\documentclass{amsart}
\usepackage{pgf}
\usepackage{amsmath}
  \usepackage{paralist}
  \usepackage{graphics} 
  \usepackage{epsfig} 
\usepackage{graphicx}  \usepackage{epstopdf}
 \usepackage[colorlinks=true]{hyperref}
\hypersetup{urlcolor=blue, citecolor=red}

\theoremstyle{definition}

\newcommand{\Rz}{\mathbb{R}}

\newcommand{\epsi}{\varepsilon} 
\renewcommand{\det}{\text{\rm det}\,} 

\newcommand{\dx}{{\rm d} x}
\newcommand{\dX}{{\rm d} X}

\newcommand{\diam}{{\rm diam}}
\newcommand{\cof}{\text{cof}\,}
\newcommand{\ove}{\overline}
\newcommand{\UUU}{\color{black}}
\newcommand{\VVV}{\color{black}}
\newcommand{\BBB}{\color{black}}
\newcommand{\RRR}{\color{black}}
\newcommand{\EEE}{\color{black}}

\newcommand{\lan}{\langle}
\newcommand{\ran}{\rangle}
\newcommand{\diameter}{\diam}

\title[Equilibrium of  immersed hyperelastic solids] 
      {Equilibrium of  immersed hyperelastic solids}

\author[Manuel Friedrich, Martin Kru\v z\'\i k, and Ulisse Stefanelli]{}


 \email{manuel.friedrich@uni-muenster.de}
 \email{martin.kruzik@utia.cas.cz}
 \email{ulisse.stefanelli@univie.ac.at}

\thanks{$^*$ Corresponding author: Ulisse Stefanelli}

\begin{document}
\maketitle

 
\medskip

\centerline{\scshape  Manuel Friedrich}
\medskip
{\footnotesize
 \centerline{Applied Mathematics,  
University of M\"{u}nster,}
\centerline{Einsteinstr. 62, D-48149 M\"{u}nster, Germany}
}

\bigskip

\centerline{\scshape Martin Kru\v z\'\i k}
\medskip
{\footnotesize
 \centerline{Academy
   of Sciences of the Czech Republic, Institute of Information Theory and Automation}
\centerline{
Pod vod\'{a}renskou v\v{e}\v{z}\'{\i}~4, CZ-182 00 Praha 8,
Czechia and}
\centerline{
Faculty of Civil Engineering, Czech Technical University,}
\centerline{Th\'{a}kurova 7, CZ--166~29 Praha~6, Czechia.}}

\bigskip

\centerline{\scshape  Ulisse Stefanelli$^*$}
\medskip
{\footnotesize
 \centerline{Faculty of Mathematics, University of Vienna,}
   \centerline{Oskar-Morgenstern-Platz 1, 1090 Wien, Austria,}
\centerline{Vienna Research Platform on Accelerating
  Photoreaction Discovery,}
\centerline{University of Vienna, W\"ahringerstra\ss e 17, 1090 Wien, Austria, and}
 \centerline{Istituto di Matematica Applicata e Tecnologie
   Informatiche {\it E. Magenes} - CNR} 
 \centerline{via Ferrata 1, 27100 Pavia, Italy.}
}
 
\bigskip


\bigskip
\bigskip

\begin{abstract}
 We discuss different equilibrium problems for hyperelastic solids 
 immersed in a fluid at rest. In particular, solids are subjected to
 gravity and hydrostatic pressure on their immersed boundaries. By
 means of a variational approach, we discuss free-floating bodies, anchored
 solids, and floating vessels. Conditions for the existence of local and
 global energy minimizers are presented.
\end{abstract}

\section{Introduction}
The equilibrium of partially immersed bodies is a classical problem in
\BBB mechanics \EEE and has attracted very early attention. Indeed, the basic
observation in this context has to be traced back more than two millennia
to the work of Archimedes  \cite{Archi}. In the first of his two books
{\it On Floating
  Bodies}, he formulated his celebrated buoyancy principle  which  is
regarded as the germinal moment of hydrostatics. In his second book, he
discusses the floating of a rigid convex paraboloid with horizontal
basis, probably inspired by the study of floating vessels.

Strangely enough, the mathematical literature on floating bodies is rather
scant. After Archimedes, an early discussion  
dates to Laplace \cite{Laplace}, who considered the case of a drop of
mercury floating on water. The rigid body case has been investigated in
\cite{John}, both at the equilibrium level and for harmonic motions. In
more recent years, a question by S.~Ulam \cite[Problem 19]{Mauldin} triggered investigations on
 the stability of convex bodies of given density
 \cite{Guerrero,Kurusa,Wegner}. Moreover, attention has been given to the capillary
case, where the fluid surface is not assumed to be flat and contact
conditions arise
\cite{Finn0,McCuan,McCuan2,McCuan3,Treinen}.  Criteria for the  stable floating of a convex
rigid body in two and
three dimensions have been recently analyzed in \cite{Finn,Finn2}. 

Driven by its obvious practical relevance,
the case of floating {\it deformable} bodies has attracted huge attention from
the engineering community. Correspondingly, the literature on {\it
  hydroelasticity} \cite{Bishop} is rather
extended. This theme fits into the general frame of {\it fluid-structure
  interactions} and the reader is referred to the recent
\cite{Kaltenbacher,Richter} for a collection of topics and references.

To the best of our knowledge, no analysis is available for
the case of a hyperelastic body deforming under the combined effect of gravity
and fluid pressure. We would like to  fill  this gap by recording in
this note some remarks on the existence
of local and global equilibria. After having collected some basic
material in Sections \ref{sec:basic}-\ref{sec:archi}, we explore a suite of different
settings, ranging from incompressible to compressible free-floating
solids (Sections \ref{sec:inc}-\ref{sec:free}), to solids at anchor
(Section \ref{sec:anchor}), to bounded fluid reservoirs (Section
\ref{sec:reserve}), to the case of ship-like
bodies (Section \ref{sec:ship}). 

In all of our discussion we follow the variational approach, by systematically restricting
our attention to energetic  arguments  and by refraining  from considering directly  the
corresponding differential problems.   On the one hand, this   reflects our personal take, which favors the relevance of
variational theories.  This, in particular, leads to recovering  a variational
version of Archimedes' celebrated principle. On the other hand, this choice allows for an
effective and
compact tractation of many different settings, which happen to be clearly
distinguished and readily amenable  by this approach.

A caveat on presentation style: in the following, we articulate a rigorous
discussion, avoiding however the classical statement-proof
structure. We hope that the reader will enjoy our informal tone, which, in
our view, better reflects the exploratory nature of our considerations.

\section{Basic setting  and energy}\label{sec:basic}

\begin{figure}[h]
  \centering
  \pgfdeclareimage[width=115mm]{figure1}{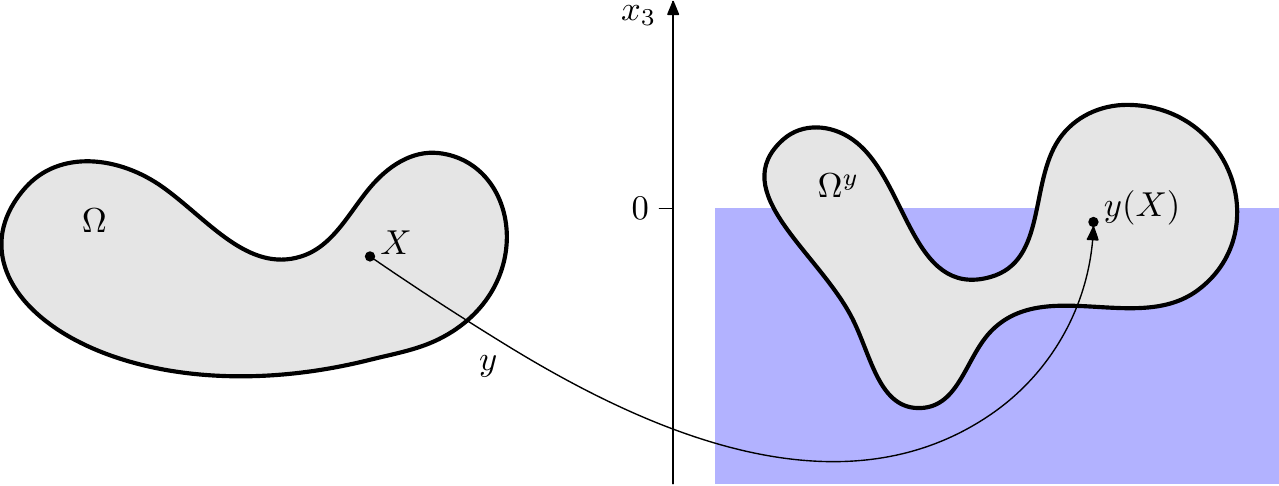}
    \pgfuseimage{figure1}
  \caption{The basic setting.}
\label{figure1}
\end{figure}

The actual shape of the solid is described by its continuous deformation $y:\Omega \to \Rz^3$ from its
reference configuration $\Omega \subset  \Rz^3$, the latter being a bounded, connected,
and smooth set. We will use
the symbol $X$ for points in the reference configuration $\Omega$ and
$x$ for points in actual space, which we endow with the orthonormal
system $(e_1,e_2,e_3)$. In coordinates, deformation $y$ is then
written as $(y_1,y_2,y_3)$.  For  a measurable set  $\omega \subset \ove
\Omega$, we denote by $\omega^y=y(\omega)$ its image under the
continuous map $y$. The actual configuration $y(\Omega)$ of the solid
is hence $\Omega^y$. 
 Both solid and fluid are  subjected  to a constant gravity in direction $-e_3$ and
 we classically indicate with $g$ the corresponding force density,
  which we assume to be constant.  In the basic setting,
 the fluid is assumed to be incompressible and to fill the
region    $\{x_3 \leq 0\} \setminus \Omega^y$, i.e.\ the part of   $\{x_3 \leq 0\}$ outside  of the body $\Omega^y$, see Figure
\ref{figure1}. 
 Correspondingly, we say that the solid {\it floats}  if $\sup y_3 >
 0$,  and that it  {\it barely floats} if $\sup y_3 = 0$. Moreover, the solid is
 said to be {\it immersed} if $\inf y_3 < 0$ and
  {\it completely immersed} if $\sup y_3 \leq 0$. 
  
The total energy of the solid is given by the sum of its {\it elastic}
potential, the {\it hydrostatic} potential, and the {\it gravitational} potential, namely 
\begin{equation}\label{eq:E}E(y) = \int_\Omega W(\nabla y(X))\, \dX + \int_{\Omega^y} g\rho_f x_3^-\, \dx +
\int_\Omega g \rho_sy_3(X)\, \dX.
\end{equation}
 We  have used the symbol $x^-=\max\{0,-x\}$ for the negative part. 
Moreover, 
we have indicated with $\rho_f>0$ the  (constant)  density of the fluid, and
with $\rho_s>0$ the  (constant)  referential density of the solid. In particular,
energy $E$ features the occurrence of both Lagrangian and Eulerian
terms. \RRR Along our discussion, we will present several variants of this basic energy  which we will highlight every time by using different notation. \EEE

Before moving on, let us justify the form of the energy 
\eqref{eq:E} by computing the Euler-Lagrange equations for $E$. Let
the invertible deformation
$y$ be given and  let  $v\colon\Omega \to \Rz^3$ represent a smooth variation.
Assuming sufficient smoothness and making use of \cite[Prop.~1.2.8,
p.~23]{Kruzik-Roubicek} for taking the variation of the hydrostatic
term we deduce that 
\begin{align*}
 \lan \delta E(y), v \ran  &= \int_\Omega {\rm D}W(\nabla y(X)):
  \nabla v(X) \, \dX \\ &+ \int_{\partial \Omega}g \rho_f y_3^-(X) \, \cof
  \nabla y(X) N(X) \cdot v(X) \, {\rm d} S(X) \\ &+\int_\Omega g\rho_s e_3 \cdot v(X)\, \dX 
\end{align*}
where  $\lan\cdot, \cdot \ran$ is  the duality product  in deformation
space,  the symbol : is the standard contraction product among
tensors,   and $N(X)$ denotes  the outward pointing normal to $\partial \Omega$.  For $y$ critical, namely  $\delta E(y) = 0$, we formally deduce the
equilibrium system for the first Piola-Kirchhoff stress $P(X)= {\rm
  D}W(\nabla y(X)) $ as
\begin{align}
  -{\rm div}_X\,P(X) +g\rho_s  e_3 = 0&\quad
  \text{in} \ \Omega,\label{eq:equilibrium}\\
P (X)N(X)= - g \rho_f y_3^-(X)\, \cof \nabla y (X) N (X)  &\quad
  \text{on} \ \partial \Omega.\label{eq:equilibrium2}
\end{align}
 Relations \eqref{eq:equilibrium}-\eqref{eq:equilibrium2} express
the equilibrium of forces in the bulk and at the boundary,
respectively. The boundary relation \eqref{eq:equilibrium2} is better understood in 
actual variables. The equilibrium system
\eqref{eq:equilibrium}-\eqref{eq:equilibrium2} can be equivalently
restated in terms of the Cauchy stress $T(y(X)) =
P(X)   (\cof \nabla y(X))^{-1}  $ as
\begin{align}
 - {\rm div}_x\,T(x) + g\rho_s  e_3  (\det \nabla
y^{-1})(x)    = 0&\quad
  \text{in} \ \Omega^y,\label{eq:equilibrium3}\\
T(x)n(x) = - g \rho_f x_3^- n(x) &\quad
  \text{on} \ \partial \Omega^y,\label{eq:equilibrium4}
\end{align}
   where  $x=y(X)$  and  $n(x)=  \cof \nabla y(X) N(X)/|\cof \nabla y(X) N(X)| $ is the outward pointing normal to $\partial \Omega^y$ at
$x$  and $y^{-1}$ is  the inverse of $y$.  Here, we used properties of Piola's transform, see \cite[Thm.~1.1.9, p.~9]{Kruzik-Roubicek}.  In particular, the hydrostatic pressure increases
linearly with depth  and no tension is exerted at its boundary 
 for $y_3\geq
0$. 

The case of a nonhomogeneous solid  could  be handled by simply letting
the density $\rho_s$ be a given integrable function  in  $X$. Moreover,
the density of the fluid can also be taken to be space-dependent. Note
that  this  does not directly apply to water, whose density varies very little
with depth, but may be relevant in  the case  of  gases. Eventually, one can
consider the case of a porous solid, which might have two densities, say $\rho_s^w$ and $\rho_s^d$,
depending on the fact that it is {\it wet} or {\it dry}. By assuming that the wet
portion of the body corresponds with the immersed one, we can model
this case by simply replacing $\rho_sy_3(X)$ by $-\rho_s^wy_3^-(X) +
\rho_s^d  y_3^+(X)  $ in the gravitational term of $E$. Here and in the
following  we let $x^+=\max\{x,0\}$. \BBB Note that here we are
excluding capillarity effects in the solid, which would wet also
portions of the solid with $y_3>0$\EEE.

\section{Requirements on the energy}
Let us specify some requirements on the elastic energy
density $W\colon{\rm GL}_+(3):=\{ F\in \Rz^{3\times 3}  \colon \, \det F >0\} \to
 [0,\infty ]  $ which will be assumed throughout. We ask $W$ \BBB to vanish \EEE at
the identity tensor $I\in
\Rz^{3\times 3}$,
\BBB and to be \EEE coercive, unbounded for $\det F\to 0+$, polyconvex, and frame indifferent. Namely, we
ask for $c_0,\, c_1>0$, $ p>3,\,s>0$, and a convex
function $\mathbb{W}:  {\rm GL}_+(3)\times  {\rm GL}_+(3) \times
  (0,\infty) \to  [0,\infty  ]  $  so that, for all $F \in {\rm
    GL}_+(3)$ we have
\begin{align}
 &   W(F)\geq c_0 |F|^p +{c_1}{(\det F) ^{-s}}- {1}/{c_0}, 
  \label{eq:coerc}\\ 
& W(F)=\mathbb{W}(F,\cof F, \det F), \\
&W(QF)=W(F) \quad \forall Q\in \text{\rm SO}(3)\label{eq:frame}
\end{align}
where $ \text{\rm SO}(3):=\{A \in \Rz^{3\times
  3}\colon\, A A^\top   = I, \ \det A=1\}$ is the set of
orientation-preserving rotations.  

Given the coercivity \eqref{eq:coerc}, finite-energy
deformations $y$ necessarily belong to the Sobolev space $ W^{1,p}(\Omega; \Rz^3)$ and
are tacitly identified with their unique continuous
representative. In particular, the set $\Omega^y$ in \eqref{eq:E} is
well-defined as the image of the
continuous function $y$. 

As the embedding of $W^{1,p}(\Omega;\Rz^3)$ in
$C(\ove \Omega;\Rz^3)$ is compact, one readily checks that energy $E$
from \eqref{eq:E} is lower semicontinuous with respect to the weak
topology of  $W^{1,p}(\Omega;\Rz^3)$. Indeed, if $y_k$  converges to $y$  weakly in
$W^{1,p}(\Omega;\Rz^3)$ and uniformly, one has that $\cof \nabla y_k
 \rightharpoonup  \cof \nabla y$ weakly in 
$L^{p/2}(\Omega;\Rz^{3\times 3})$,  $\det \nabla y_k   \rightharpoonup  \det \nabla y$ weakly in
$ L^{p/3}  (\Omega)$, and $((y_k)_3)^{-} \to y_3^-$ uniformly, so that 
\begin{align*}
  \liminf_{k \to \infty} E(y_k) &=  \liminf_{k \to \infty}
  \Bigg(\int_\Omega \mathbb{W}(\nabla y_k(X) ,\cof \nabla y_k(X) ,\det \nabla
  y_k(X) ) \, \dX \nonumber\\
&+ \int_\Omega g\rho_f y_3^-(X) \det \nabla y(X) \, \dX + \int_\Omega g \rho_s y_3(X) \, \dX
  \Bigg)\geq E(y).
\end{align*}
Having settled lower semicontinuity, the discussion on existence of
energy minimizers will focus on identifying conditions ensuring the coercivity
of the energy $E$ with respect to the weak topology of
$W^{1,p}(\Omega;\Rz^3)$. Indeed, once such coercivity  is established,
existence of global energy minimizers
would follow by the \BBB direct method\EEE. 

 
\VVV In the following, we \EEE ask admissible deformations $y\in W^{1,p}(\Omega;\Rz^3)$ to additionally
fulfill the classical Ciarlet-Ne\v cas  \RRR condition \BBB 
\cite{CN} 
\begin{equation}\label{eq:CN}
\int_\Omega \det \nabla y(X) \, \dX \leq |\Omega^y|,
\end{equation}
 where, here and below
$|\cdot|$ stands also for the Lebesgue measure in $\Rz^3$. This
condition implies that the map $y$ is almost everywhere injective,
namely, there  exists  $\omega\subset\Omega$ such that $|\omega|=0$ and 
$y(x_1)\not = y(x_2)$  for every  $x_1,x_2\in\Omega\setminus\omega$
satisfying $x_1\not = x_2$. Note that \eqref{eq:CN} is 
closed under weak $W^{1,p}(\Omega;\Rz^3)$ \RRR convergence. \BBB

As almost everywhere injective can still be noninjective, one \VVV
could consider strengthening \EEE 
the coercivity of the energy \EEE as 
$$W(F)\geq c_0 |F|^p + c_0 \frac{|F|^{3q}}{(\det
  F)^q}+{c_1}{(\det F) ^{-s}}- \frac{1}{c_0}$$
for some $q>2$. 
Under this stronger coercivity condition, finite-energy deformations are of
  {\it finite-distortion} \cite[Def.~1.11, p.~14]{HK}. In particular, owing to
\cite[Thm.~3.4, p.~43]{HK} they are {\it open}, \BBB namely, map open
sets to open sets. The last step is then to check that almost
everywhere injective open deformations $y$ are actually homeomorphisms
from $\Omega$ to $\Omega^y$ \cite[Thm.~3.5]{GKMS}. This in particular,
entails that $y$ is injective\EEE.

\section{Coercivity and the Archimedes Principle}\label{sec:archi}

As already mentioned, coercivity of the energy in the weak  $W^{1,p}(\Omega;\Rz^3)$ 
topology entails existence of minimizers.
\BBB Note that the \EEE elastic energy controls the $L^p$ norm of $\nabla y$ via
\eqref{eq:coerc}. \BBB Thus, \EEE in order to deduce the coercivity of $E$ with
respect to the weak $W^{1,p}(\Omega;\Rz^3)$ topology, we just need to ascertain that admissible deformations are uniformly bounded \BBB on \EEE some energy sublevel. The elastic part of the energy is 
invariant under rigid motions and,  since $p>3$,   an  $L^p$ bound on $\nabla y$ entails a 
bound  on the {\it diameter} $\diameter(\Omega^y)$ of $\Omega^y$, namely
\begin{equation}
\diameter(\Omega^y) \leq c \|\nabla
y\|_{L^p(\Omega)}. \label{eq:diam}
\end{equation}
Here and in the following, we will use the symbols $c$, $c'$,
$c''$ to
indicate positive constants, possibly depending on data but
independent of the deformation, and changing from line to line. 

Owing to the diameter bound \eqref{eq:diam}, the boundedness of $y$ will follow
as soon as one checks that the position of the barycenter $\ove y $ of the
\BBB deformed body \EEE $$\ove y = (\ove y_1,\ove y_2,\ove
y_3) { := } \frac{1}{|\Omega|}\int_\Omega
y(X)\, \dX \in \Rz^3$$
is bounded by the energy.  Since the hydrostatic and
the gravitational terms depend just on the
$y_3$ component,  it is not restrictive to assume that 
$$\ove y_1 = \ove y_2 = 0.$$
Hence,  one needs to check the boundedness of $\ove y_3$ only.

The gravitational potential strictly decreases
by translating the solid in direction $-e_3$.
As the hydrostatic term vanishes out of the
fluid, i.e., for $y_3>0$, one readily proves that  solids
minimizing the energy  are necessarily immersed,  namely  $\inf y_3 <0$. \BBB This in turn implies $\sup y_3 \le \inf y_3 + \diameter(\Omega^y) < \diameter(\Omega^y)$.  By definition of the barycenter this yields \EEE  
\begin{equation}
\ove y_3
<\diameter(\Omega^y).\label{eq:upper}
\end{equation} 
The issue is then to control $\ove y_3$ from below. Assume then that
the  solid   is completely immersed and rewrite $E$ as
\begin{align}
  E(y) &= \int_\Omega W(\nabla y(X))\, \dX + g \int_\Omega \left(-\rho_f
 J^y(X) + \rho_s\right) y_3(X)\, \dX\nonumber\\
&= \int_\Omega W(\nabla y(X))\, \dX + g \int_\Omega \left(-\rho_f
 J^y(X) + \rho_s\right) (y_3(X)-\ove y_3)\, \dX\nonumber \\
&\ \ \ + g |\Omega| \left(-\rho_f
 \ove J^y + \rho_s\right) \ove y_3,\label{eq:Eimm}
\end{align}
where we have used the short-hand notations
$$J^y(X) = \det \nabla y(X) \quad \text{and} \quad \ove J^y =
\frac{1}{|\Omega|} \int_\Omega J^y(X)\, \dX.$$
The first two terms in the right-hand side of \eqref{eq:Eimm} are
invariant under translations in direction $e_3$.   The last term  is decreasing  as $\ove
y_3$ increases iff 
\begin{equation}
  \label{eq:buoy}
  \rho_f
 \ove J^y> \rho_s.
\end{equation}
 Relation \eqref{eq:buoy}  is hence a necessary and sufficient buoyancy
  condition.

Let now $y$ be a critical point for $E$ and consider variations 
of the form 
$  y+\alpha e_3$ for $\alpha \in \Rz$. From $f'(0) =0 $ for
$f(\alpha) = E(y+\alpha e_3)$ we deduce  from \eqref{eq:E}  that 
\begin{equation}
 g \rho_f \, |\Omega^y\cap \{x_3 \leq 0\}| =   g
\rho_s \, |\Omega|.\label{eq:archi}
\end{equation}
The left and the right term in this equation are respectively the weight of
 the displaced fluid and the weight of the solid, so that relation
 \eqref{eq:archi} is nothing but the classical {\it Archimedes
principle}. 

\section{Incompressible solids}\label{sec:inc}

 We consider an {\it
  incompressible} free-floating solid  by requiring $J^y=1$ almost
everywhere  (for instance by
letting $W(F)=\infty$ if $\det F \not =1$).  Note that the
incompressibility constraint  $J^y=1$ a.e.\ is stable  under  weak $W^{1,p}(\Omega;\Rz^3)$ 
convergence.  In this case,  condition   \eqref{eq:buoy}
reduces to  a relation between  $\rho_s$ and $\rho_f$.  

 In case $\rho_f> \rho_s$ we have that  solids
minimizing the energy  
necessarily float. In
particular, the energy controls the full $W^{1,p}$ norm of the
deformation and the \BBB direct method \EEE ensures the existence of a ground
state $y$. One can  check  \BBB that \EEE  $\Omega^y$  is  floating,  for one has  
$$|\Omega^y \cap \{x_3\leq 0\}| \stackrel{\eqref{eq:archi}}{=} (\rho_s/\rho_f)|\Omega|<
|\Omega|.$$

If $\rho_f<\rho_s$, the energy $E$ is not bounded from below since it
 decreases linearly for 
translations of the solid in direction $-e_3$  once the  solid  is completely immersed,  see \eqref{eq:Eimm}.  In this case, no global minimizer of
$E$ exists. 

In the critical case $\rho_f=\rho_s$, one can rewrite the energy as
$$E(y) = \int_\Omega W(\nabla y(X))\, dX + \int_\Omega g\rho_s
y^+_3(X)\, \dX \geq 0.$$ 
In particular, completely immersed, 
rigid  solids,  namely  those given by  $y(X) = QX + v$ with $Q \in {\rm SO}(3)$,
$v \in \Rz^3$, and $\sup y_3 \leq 0$,
realize $E=0$ and are thus global minimizers. In case $W=0$ solely on
${\rm SO}(3)$, these are actually the unique
global minimizers.

The existence of   floating  or barely floating  solids minimizing
the energy 
can be proved even in the case of  slight compressibility,  as
long as $\ove J^y$ is bounded away
from zero, say \begin{equation}
  \label{eq:buoy2}
  \ove J^y\geq  \tau  >0.
\end{equation} In this case, one can still recover the existence 
of  floating  solids minimizing
the energy 
 whenever  $
  \rho_f \tau \geq \rho_s$.
A lower bound on $\ove J^y$ of the form of \eqref{eq:buoy2} would follow
in case $W$  was  constrained to be $W(F)=\infty$ for $\det F < \tau$
(a closed condition with respect to the weak $W^{1,p}$
topology). Alternatively, one could replace the term $c_1(J^y)^{-s}$ by
$c_1(J^y-\tau)^{-s}$ in the coercivity condition \eqref{eq:coerc} in
order to 
restrict to deformations with $J^y>\tau$  almost everywhere. 

The findings of this section  for incompressible solids  can be summarized as follows:
\begin{itemize}
\item If $\rho_s > \rho_f$ no local minimizer exists and the energy
  is not bounded from below.
\item If $\rho_s = \rho_f$  solids minimizing
the energy  exist and are 
  completely immersed.
\item If $\rho_s < \rho_f$  solids minimizing
the energy  exist and are floating. 
\end{itemize}

\begin{figure}[h]
  \centering
  \pgfdeclareimage[width=55mm]{submarine}{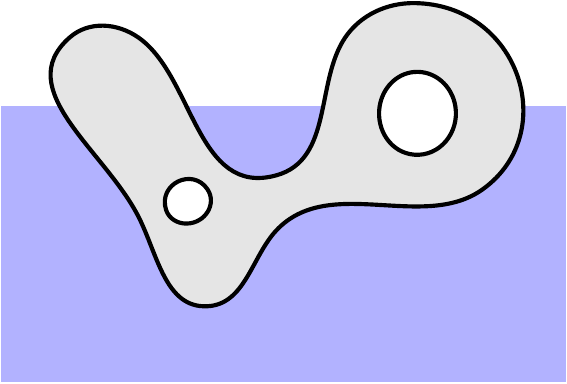}
    \pgfuseimage{submarine}
  \caption{The submarine setting.}
\label{submarine}
\end{figure}

By possibly resorting to a nonhomogeneous density $\rho_s \colon\Omega \to
[0,\infty)$ one could discuss the case of  a  vessel having load or
flotation tanks. This is for instance the case of a submarine, see
Figure \ref{submarine}, where buoyancy is controlled by allowing water to
fill the ballast tanks or by expelling water from the ballast tanks by
means of a compressed-air reserve.
 In order to model the floating of a submarine, one assumes $\Omega$ to be disjointly partitioned as $ \Omega=  \Omega_h
 \cup  \Omega_b $, where  $\Omega_h$  is the reference configuration of
the hull of the submarine (where density  $\rho_s=\rho_h$  is assumed to be
constant) and $\Omega_b$ represents the ballast tanks  with  density $\rho_s=\rho_b$ depending on the air-water ratio.  A possible form
for the energy in this case is
\begin{align}
\RRR E_{\rm submarine}(y) \EEE &= \int_{\Omega_h}W(\nabla y(X))\, \dX +
\int_{\Omega^y} g \rho_f x_3^-\,\dx\nonumber
\\
&\quad+\int_{\Omega_h}g \rho_h y_3(X)\, \dX +\int_{\Omega_b}g \rho_b y_3(X)\, \dX.\label{eq:E6}
\end{align}
 In order to
determine the correct air-water balance in the ballast tanks keeping
the  (incompressible)  submarine completely immersed and
neutrally floating,  one has to simply 
reconsider the discussion of \eqref{eq:archi} 
to find the density $\rho_b$ (assumed constant for simplicity) fulfilling
$\rho_f |\Omega|  = \rho_h |\Omega_h| + \rho_b
|\Omega_b|$.  (Here, neutrally floating means that the energy is independent from $\ove
y_3$ as long as the  solid  is fully immersed.)

\section{Compressible  solids}\label{sec:free}
Let us now turn to the case of a free-floating compressible solid. In this case, regardless of the values of
$\rho_f$ and $ \rho_s$, the energy can be proved to
be unbounded from below, so that no global minimizer exists. Indeed, for
all $\rho_f,\, \rho_s>0$, one can consider a deformation $y$ with
$\rho_f
        J^y <\rho_s$ almost everywhere, define $y_k = y - k e_3$ for $k \geq
        \diameter(\Omega^y)$, so that  $\Omega^{y_k}$  is completely immersed, and
compute 
\begin{align}
 E(y_k) &=\int_\Omega W(\nabla y(X))\, \dX + \int_\Omega g(-\rho_f
        J^y(X) + \rho_s)(y_3(X)-k)\, \dX \nonumber \\
&\leq \int_\Omega W(\nabla y(X))\, \dX+ g(-\rho_f
        \ove J^y + \rho_s) |\Omega| ( \sup  y_3 -k).\label{eq:ab}
\end{align}
By taking $k \to \infty$ one has that $E(y_k) \to -\infty$. Note that
the energy is unbounded  from  below for any choice of the densities,
regardless of the fact that $\rho_s$ could be smaller than
$\rho_f$ (and hence the solid would float,  were it
incompressible). 

As global minimizers do not exist, we now turn to consider local
minimizers instead. In case $\rho_s>\rho_f$, given any $y$  with $J^y = 1$ almost everywhere,  one has
that $\alpha \mapsto E(y+\alpha e_3)$ increases  on  $\lbrace \alpha<0\rbrace$.  In particular, $E$ admits no local
energy minimizer.

We hence turn to the case $\rho_s < \rho_f$ for the remainder of
this section.
The argument in \eqref{eq:ab} is based on considering
deformations with
$J^y(X)<\rho_s/\rho_f$. In fact, $J^y(X)$ can be made arbitrarily
small, still keeping the energy finite. On the other hand, the
coercivity \eqref{eq:coerc} entails that extreme compressions have high
energy. Indeed, the higher the value of the parameter
$c_1$, the higher is the energy needed to obtain
$J^y(X)<\rho_s/\rho_f$. It is hence conceivable that, for \BBB
materials having a \EEE very large $c_1$
and for $\rho_s<\rho_f$, local minimizer  of the energy  among
floating  solids  may exist. We devote the following discussion to check this
fact,  by assuming a  slightly more specific form of the
energy,  namely, 
\begin{align}
\RRR E_{\rm sp}(y) \EEE &= \int_\Omega W(\nabla y(X)) \, \dX + c_1 \int_\Omega
(J^y(X))^{-s}\, \dX \nonumber\\
&\quad+ \int_{\Omega^y} g\rho_f x_3^-\, \dx +
\int_\Omega g \rho_sy_3(X)\, \dX\label{eq:E4}
\end{align}
for $W\colon{\rm GL}_+(3) \to [0,\infty]$ with $W(I)=0$  fulfilling
\eqref{eq:coerc}-\eqref{eq:frame}  (with $c_1=0$ \BBB there). \EEE In particular, we
have  highlighted  the coercive part on $(J^y)^{-s}$  by
separating it  from the elastic
energy.

Before moving on, we need to refine the bound \eqref{eq:diam} for a
floating solid,  i.e., for  deformations  with  $\sup y_3 >0$. \BBB By   \eqref{eq:coerc}  (for $c_1=0$) \EEE we have that 
\begin{align*}
\diameter(\Omega^y) & \leq c \| \nabla y\|_{L^p(\Omega)} \leq c'\left( \RRR E_{\rm sp}(y) \EEE + 1/c_0 -
  \int_{\Omega^y}g \rho_f x_3^-\, \dx - \int_\Omega g\rho_s y_3(X)\,
  \dX \right)^{1/p}\\
&\leq c'\big( \RRR E_{\rm sp}(y) \EEE +  1/c_0  + c'' \diameter(\Omega^y) \big)^{1/p}.
\end{align*}
From this we get that 
\begin{equation}
  \label{eq:diam2}
\diameter(\Omega^y) \leq c\left( \RRR E_{\rm sp}(y) \EEE \right)^{1/p} +c.
\end{equation}

Let us start by restricting ourselves to a specific sublevel of the
energy. Define $\hat y:= {\rm id} + \alpha e_3$, where $\alpha\in \Rz$ is
chosen in such a way that $\sup \hat y_3=0$. \BBB In this case, it holds that  $W(\nabla \hat{y}(X)) =0$ and $J^{\hat{y}}(X)=1$ for all $X \in \Omega$ by $W(I)=0$. Moreover, we have $x_3^- = -x_3$ for all $x \in \Omega^{\hat y}$. Then, \EEE we readily compute that
\begin{equation}
  \label{eq:Ebar}
  \RRR E_{\rm sp}(\hat y) \EEE = c_1 |\Omega| + g(\rho_s - \rho_f)(\ove X_3 +\alpha)|\Omega|
\end{equation}
where  $\ove X_3:= \int_\Omega X_3 \dX /|\Omega|$.   As $\rho_s<\rho_f$,
such $\hat y$ cannot be a minimizer of the energy, for the energy
decreases by increasing $\alpha$.  As a first step,  we fix $1> \tau> \rho_s/\rho_f$ and we   check that the energy $ \RRR E_{\rm sp}\EEE$ can be minimized on the
set
$$A=\{ y \in W^{1,p}(\Omega;\Rz^3) \ \ \text{fulfilling  $\ove
  J^y\geq \tau$}\}.$$
In fact,  for  minimizers we can restrict
our considerations to the sublevel $\{ \RRR E_{\rm sp} \EEE \leq   E_{\rm sp}(\hat y) \EEE\}$. All such
 solids  are necessarily floating, see \eqref{eq:buoy2} in  Section
\ref{sec:inc}.  Thus, in particular, the energy $ \RRR E_{\rm sp}\EEE$ is coercive on  $A \cap \{ \RRR E_{\rm sp} \EEE \leq  \RRR E_{\rm sp}(\hat y) \EEE\}$  and
$A$ is closed with 
respect to the weak $W^{1,p}$-topology.  This implies that there
exists a  
minimizer $y^*\in A$ of $ \RRR E_{\rm sp}\EEE$ on $A$  (with $\Omega^{y^*}$
floating).

 Our key step is to  check that, \BBB if $c_1>0$ is sufficiently
 large, \EEE 
one can find $r>0$ small such that 
\begin{equation}
  \label{eq:pro}
  \big(\| y_3 - y^*_3 \|_{L^\infty (\Omega) }<r, \  \RRR E_{\rm sp}(y)\leq  E_{\rm sp}(\hat y)\EEE\big) \ \Rightarrow \
  \ove J^y >  \tau. 
\end{equation}
 This then  proves that $y^*$ is a local  (in $L^\infty$)  
minimizer of the energy.  In fact, given $y$ with $\| y_3 - y^*_3 \|_{L^\infty(\Omega)}<r$, one either has $ \RRR E_{\rm sp}(y) >   E_{\rm sp}(\hat{y}) \ge  E_{\rm sp}(y^*)\EEE$ or $ \RRR E_{\rm sp}(y)\leq  E_{\rm sp}(\hat y)\EEE$. In the latter case,   implication \eqref{eq:pro} entails that $y \in A$ and therefore $ \RRR E_{\rm sp}(y) \ge E_{\rm sp}(y^*)\EEE$.

In order to check \eqref{eq:pro}, we argue by contradiction and assume
to be given a $ y \in W^{1,p}(\Omega;\Rz^3)$ with  $\| y_3 - y^*_3 \|_{L^\infty (\Omega) }<r$ and   $ \RRR E_{\rm sp}(y)\leq   E_{\rm sp}(\hat y)\EEE$ such that 
\begin{equation}
  \label{eq:tilde}
 \ove J^y\leq  \tau. 
\end{equation}
By letting $ \mu  =g(\rho_s - \rho_f)(\ove X_3 +\alpha)$  and recalling that $W \ge 0$,  we then have
\begin{align*}
   &c_1 |\Omega| + \mu |\Omega|
     \stackrel{\eqref{eq:Ebar}}{=}  \RRR E_{\rm sp}(\hat y) \geq  E_{\rm sp}(y)\EEE \\
&\quad  \stackrel{\eqref{eq:E4}}{\geq}
     c_1\int_\Omega (J^y(X))^{-s}\, \dX + \int_{\Omega^y} g \rho_f x_3^-
      + \int_\Omega g \rho_s y_3(X)\, \dX\\
&\quad \geq c_1|\Omega| \big(\ove J^y\big)^{-s}  + \int_{\Omega^{y^*}} g \rho_f x_3^-
     + \int_\Omega g \rho_s y^*_3(X)\, \dX - c\| y_3 - y^*_3 \|_{L^\infty (\Omega)  }\\
&\quad \stackrel{\eqref{eq:tilde}}{\geq} c_1|\Omega|\tau^{-s}- g
  \rho_s|\Omega|\diameter(\Omega^{y^*}) - cr\\
&\quad
  \stackrel{\eqref{eq:diam2}}{\geq}c_1|\Omega|\tau^{-s}
  -c ( \RRR E_{\rm sp}(y^*)\EEE)^{1/p} - c - cr\\
&\quad \geq c_1|\Omega|\tau^{-s} -c ( \RRR E_{\rm sp}(\hat y)\EEE)^{1/p} -
  c  - c r\\
&\quad \stackrel{\eqref{eq:Ebar}}{=}
  c_1|\Omega|\tau^{-s} -c (c_1 |\Omega| +  \mu  |\Omega|)^{1/p} -
  c - c r.
\end{align*}
 Here, we have also used Jensen's inequality in the third inequality  and $\sup y_3^* >0$ in the fourth inequality. 
In the latter computation and up to the end of this section, the
generic  constant $c$ is  \BBB always \EEE independent of $c_1$ and
$r$ as well. We have checked that 
\begin{equation}c_1+  \mu  \geq c_1 \tau^{-s} -c c_1^{1/p}-c - c r.\label{eq:contra}
\end{equation}
 Since $\tau<1$,  given any $r >0$ the latter does not hold \BBB in
 case $c_1>0$ is sufficiently \EEE large. 
This \BBB leads \EEE to a contradiction, proving
\eqref{eq:pro}. 

In conclusion, we have checked that, for all $\rho_s<\rho_f$ and all
$r>0$, there exists a $c_1>0$ such that a minimizer of $ \RRR E_{\rm sp} \EEE$ in $A$
is a local minimizer of $ \RRR E_{\rm sp} \EEE$ (in
a $L^\infty$ ball of radius $r$).

 Note that,  given any $c_1$, the argument
of \eqref{eq:contra} fails for $r$ large enough.  This corresponds to
the former observation that the energy is de facto unbounded from below. On the other hand, for small values of $c_1$, \eqref{eq:contra} does not allow to
conclude for the existence of $r$ such that a local minimizer in
the $L^\infty$ ball of radius $r$ exists. In fact, in the limiting
case $c_1=0$  and    $W=0$, one can check that the
energy $ \RRR E_{\rm sp} \EEE$ has not even local minimizers. 

The discussion on the buoyancy of a submarine from Section \ref{sec:inc}, see Figure
\ref{submarine}, can be extended to the case where the air-water  balance  of the
ballast tanks is compressible (but, for simplicity, the hull of the
submarine is not). In this case, we could specify the energy from
\eqref{eq:E4} by following \eqref{eq:E6}, namely,
\begin{align*}
\RRR E_{\rm submarine, sp}(y) \EEE &= \int_{ \Omega_h } W(\nabla y(X))\,
  \dX  + c_1\int_{\Omega_b}(J(X))^{-s}\dX + \int_{\Omega^y} g\rho_f x_3^-\,
                     \dx \nonumber\\
&\quad +
\int_{\Omega_h} g \rho_hy_3(X)\, \dX + \int_{\Omega_b} g \rho_by_3(X)\, \dX.
\end{align*}
The argument above can be adjusted to the case of the latter
energy. \BBB In case $c_1$ and $|\Omega_b|$ are sufficiently large, \EEE one can find  neutrally floating 
 solids locally minimizing the energy.

\section{Solids at anchor}\label{sec:anchor}
\begin{figure}[h]
  \centering
  \pgfdeclareimage[width=115mm]{anchoring}{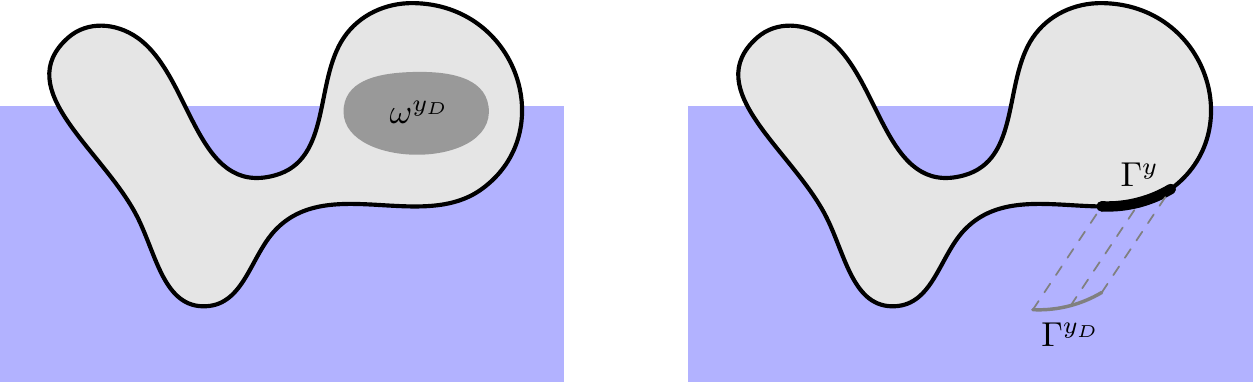}
    \pgfuseimage{anchoring}
  \caption{Two anchored situations: prescribed deformation on $\omega
    \subset \Omega$ (left) and elastic boundary conditions on
    $\Gamma\subset \partial \Omega$ (right).}
\label{anchoring}
\end{figure}

Independently of compressibility,  the  coercivity of the energy in
$L^\infty$ (and hence in the weak topology of $W^{1,p}$) can be obtained by
prescribing some form of {\it anchoring} of the solid. A classical choice in
this sense consists in assuming a prescribed deformation $y = y_{D}$ at some point in $\Omega$ (or on some portion $\omega$ of
$\Omega$ or  some portion $\Gamma$ of  $\partial \Omega$), see Figure \ref{anchoring} left. 
In this case, one can compute
\begin{align}
  |\Omega|/c_0  +  E(y) &\geq c_0 \| \nabla y \|^p_{L^p(\Omega)}  
+ \int_{\Omega^y} g\rho_f x_3^- \, \dx + \int_\Omega g\rho_s y_3(X)
                  \, \dX\nonumber \\
& \geq c_0 \| \nabla y \|^p_{L^p(\Omega)}    -  g \rho_s(\diameter
  (\Omega^y) +|y_{D}|) \nonumber\\
&\stackrel{\eqref{eq:diam}}{\geq} c_0 \| \nabla y \|^p_{L^p(\Omega)}    - c
  \|\nabla y\|_{L^p(\Omega)}   - c \geq \frac{c_0}{2} \| \nabla y \|^p_{L^p(\Omega)}    -
   c' . \label{eq:previous}
\end{align}
In particular, this entails that 
$${ \| y_3 \|_{L^\infty(\Omega)}   \leq \diameter(\Omega^y) + |y_{D}|
\stackrel{\eqref{eq:diam}}{\leq} c \| \nabla y\|_{L^p(\Omega)}  + |y_{D}|
\stackrel{\eqref{eq:previous}}{\leq} c(E(y))^{1/p} +  c'  + |y_{D}| }$$
so that the energy is coercive in $L^\infty$. However, one has to mention that prescribing some specific deformation
at some point or portion of $\Omega$ or $\partial \Omega$ could be
practically not realizable, especially in three space dimensions, see
again Figure \ref{anchoring} left.

An alternative choice could be that of considering so-called {\it elastic} boundary conditions,
again to be imposed  on \BBB a portion of $\partial\Omega$ (or,
alternatively, on a portion of $ \Omega$). \EEE In order to give an example in this
direction, we fix $\Gamma \subset \partial \Omega$ open in the
relative topology of $\partial \Omega$ with $\mathcal H^2(\Gamma)>0$,
let $y_{D}\colon \Gamma \to \Rz^3$ be continuous, and augment the energy
by the term
\begin{equation} 
c_3\int_\Gamma  |y(X) - y_{D}(X)|^r\, {\rm d} \mathcal
H^2(X),\label{eq:boundary}
\end{equation}
\RRR denoted by $E_\Gamma$. \EEE  Here,  $\mathcal H^2$ denotes the two-dimensional Hausdorff measure on
$\Gamma$, $c_3>0$, and $r> 1$. See Figure \ref{anchoring} right. By taking into account the Poincar\'e inequality  (see, e.g.,
\cite[Lemma 3.3]{compos2}),  we obtain
\begin{align}\label{poinci}
\|y\|_{W^{1,p}(\Omega)  } \leq c \|\nabla
y\|_{L^p(\Omega)} + c\|y - {\rm id}\|_{L^1(\Gamma)} \quad \forall \,  y \in
W^{1,p}(\Omega;\Rz^3).  
\end{align}
\BBB We can adapt the argument in
\eqref{eq:previous} to this situation:  first, we get \EEE
\begin{align*}
 |\Omega| /c_0  +  \RRR E_\Gamma(y) \EEE &\geq c_0 \| \nabla y \|^p_{L^p(\Omega)}  +c_3\| y - y_{D}\|_{L^r(\Gamma)}^r
+ \int_{\Omega^y} g\rho_f x_3^- \, \dx + \int_\Omega g\rho_s y_3(X)
                  \, \dX\nonumber \\
& \geq c_0 \| \nabla y \|^p_{L^p(\Omega)}  +c_3\| y -
  \BBB {\rm id} \EEE \|_{L^r(\Gamma)}^r - \BBB c'' \EEE - 
  g\rho_s  \|y\|_{L^\infty(\Omega)}.
\end{align*}
\BBB Then, due to  \eqref{poinci}  and the embedding $W^{1,p}\subset L^\infty$, we obtain \EEE
\begin{align}
 |\Omega| /c_0  +   \RRR E_\Gamma(y) \EEE & \geq c \|  y \|^{\min\{p,r\}}_{W^{1,p}(\Omega)}  - c' \|  y \|_{W^{1,p}(\Omega)}
  -c''\nonumber\\
&\geq c \|  y \|^{\min\{p,r\}}_{W^{1,p}(\Omega)} -c,\label{eq:previous2}
\end{align}
\BBB where in the last step we used that  $c z^q -c' z \ge - c_*$  for all $z \in (0,\infty)$, provided $q>1$ and $c_*$ is sufficiently large. \EEE Once again, the energy is coercive in $L^\infty$.

A further sophistication could be that of assuming that the elastic
response of the boundary condition is inactive before a given
critical elongation $|y - y_{D}|=\lambda>0$ is reached. This would indeed correspond
to the case of a buoy or a floating vessel anchored via an elastic
cable  of length $\lambda$ at rest.  In this case, the boundary condition term would be modified as 
$$ c_3\int_\Gamma  \max\{0,|y(X) - y_{D}(X)|^r-\lambda^r\} \, {\rm d} \mathcal H^2(X).$$
The argument in \eqref{eq:previous2} can be adapted to this case by simply
replacing the term $c_3\| y - y_{D}\|^r_{L^r(\Gamma)}$ by $c_3\| y -
y_{D}\|^r_{L^r(\Gamma)} - c_3\mathcal H^2(\Gamma)\lambda^r$.

Eventually, one could consider the case of an inextensible
anchoring of length $\lambda>0$. This would be modeled by imposing the constraint 
$$ \| y - y_{D} \|_{L^\infty} \leq \lambda$$
which would directly entail coercivity, simply by replacing $|y_{D}|$ by
$|y_{D}|+\lambda$ in the chain of inequalities \eqref{eq:previous}.

\section{Bounded reservoir}\label{sec:reserve}
\begin{figure}[h]
  \centering
  \pgfdeclareimage[width=60mm]{reserve}{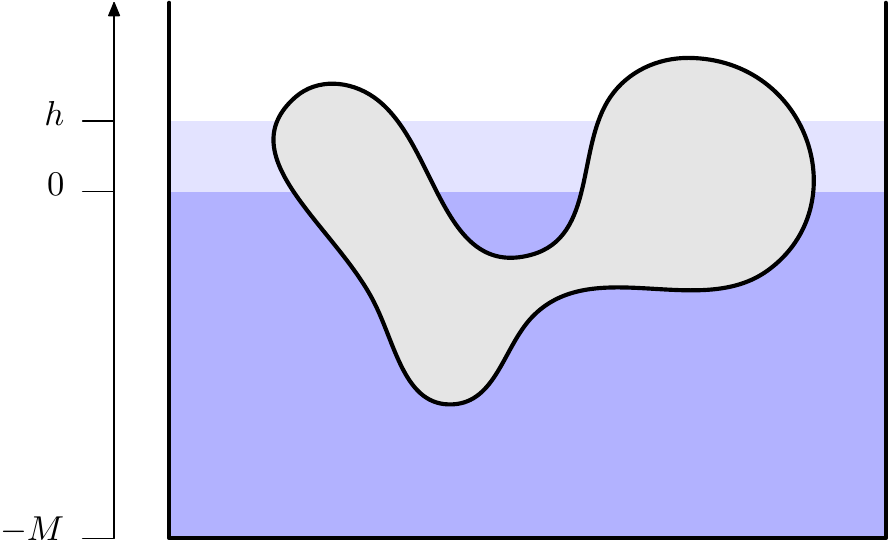}
    \pgfuseimage{reserve}
  \caption{The bounded-reservoir setting.}
\label{reserve}
\end{figure}

In this section we discuss the case of a \UUU compressible \EEE solid floating 
in a {\it bounded} fluid reservoir. We assume the container to be large
enough so that the solid does not touch its walls, see Figure
\ref{reserve}. 
For the sake of definiteness, we assume the reservoir to
 be of cylindrical
shape $S\times [-M,\infty)$ with $S\subset \Rz^2$ compact and $M>0$ sufficiently large. 
 The issue is here that the immersed portion of the
solid lets the  water level  rise to some level $h$, so that the energy
$E$ from \eqref{eq:E} has to be modified as 
\begin{equation}\label{eq:E3}
\RRR E_{\rm reservoir}(y,h) \EEE = \int_\Omega W(\nabla y(X))\, \dX + \int_{\Omega^y} g\rho_f (x_3-h)^-\, \dx +
\int_\Omega g \rho_sy_3(X)\, \dX.
\end{equation}
Note that the energy now depends on $h$ as well, and this is a priori
unknown.

In case the  solid  is free-floating, see Sections
\ref{sec:inc}-\ref{sec:free}, the treatment of the extra variable $h$
is straightforward: letting  $y_*$  be a (either
local or global)  minimizer of the original energy from \eqref{eq:E},
one has that  $y:=y_*+he_3$  minimizes $\RRR E_{\rm reservoir} \EEE$ from \eqref{eq:E3}, where $h$ is a
posteriori determined by solving
\begin{equation}\label{eq:cont}
M \mathcal H^2(S) +  |\Omega^{y}\cap \{x_3 \leq h\}|  =(M+h)\mathcal H^2(S).
\end{equation} 
This equation is nothing but the expression of the conservation of fluid
content. On the left is the sum of the original and the displaced
fluid volume. On the right is the volume of the portion of the
reservoir under water level.

The problem is more involved if the solid is anchored, see
Section \ref{sec:anchor}. In this case, the value of $h$ cannot be
computed a posteriori and one has to minimize $\RRR E_{\rm reservoir} \EEE$ from
\eqref{eq:E3} directly,   under the additional constraint \eqref{eq:cont}. 
This is 
however possible, as we now check. Assume for definiteness that the
solid is clamped, as in the left of Figure \ref{anchoring} (other
cases, including \eqref{eq:boundary}, can be treated as well). Let $(y_k,h_k)$ be an infimizing
sequence for
\eqref{eq:E3} under the constraint \eqref{eq:cont}. The values $h_k$
are surely bounded from below by $0$, for the water level can only raise as
effect of the immersed  solid.   In a similar fashion, the values  $h_k$ are bounded from
above by $|\Omega|\ove J^{y_k}/ \mathcal{H}^2(S)$, which in turn is bounded by the
energy itself. Up to  subsequences (not relabeled), we hence have that  
$y_k \rightharpoonup y$ weakly in $W^{1,p}(\Omega;\Rz^3)$ and $h_k \to h$ in $\Rz$.  In particular, this entails that   $
g\rho_f(x_3-h_k)^- \to  g\rho_f(x_3-h)^-$ uniformly
and $|\Omega^{y_k} \triangle \Omega^y| \to 0 $ \cite[Lemma 5.2]{GKMS},  where $\triangle$ denotes the symmetric difference of sets.  The
latter implies that $
\chi_{\Omega^{y_k}}\to \chi_{\Omega^{y}}$ strongly in $L^1$
where $\chi_{\RRR A}(x)$ is the characteristic function of the measurable set
$\RRR A \EEE \subset \Rz^3$. We can hence
pass to the limit in the hydrostatic terms and obtain 
\begin{align*}
  &\int_{\Omega^{y_k}}g \rho_f (x_3-h_k)^- \, \dx = \int_{\Rz^3}g
  \rho_f (x_3-h_k)^- \, \chi_{\Omega^{y_k}}(x)\,\dx\\
&\quad \to  \int_{\Rz^3}g
  \rho_f (x_3-h)^- \, \chi_{\Omega^{y}}(x)\,\dx= \int_{\Omega^{y}}g \rho_f (x_3-h)^- \, \dx.
\end{align*}
This convergence entails the lower semicontinuity
of the energy $\RRR  E_{\rm reservoir} \EEE$. In addition, 
$|\Omega^{y_k} \triangle \Omega^y| \to 0 $  and $h_k \to h$ entail  that 
$$ |\Omega^{y_k} \cap \{x_3 \leq  h_k  \}| \to |\Omega^{y} \cap \{x_3 \leq
 h  \}|.$$ 
One can hence pass to the limit in equation \eqref{eq:cont}, written
for $y_k$, in order to check that the limiting pair $(y,h)$ fulfills
\eqref{eq:cont} as well.

\section{The ship problem}\label{sec:ship}
\begin{figure}[h]
  \centering
  \pgfdeclareimage[width=55mm]{ship}{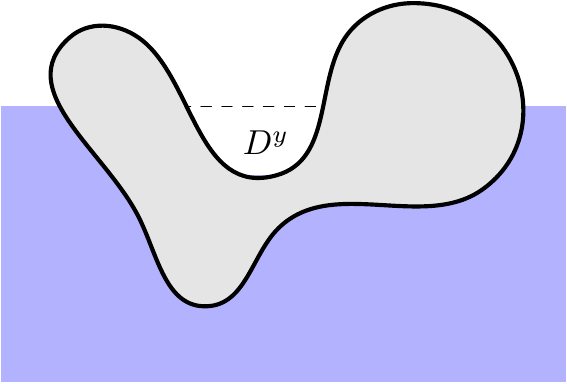}
    \pgfuseimage{ship}
  \caption{The ship setting.}
\label{ship}
\end{figure}

Let us now go back to the case of an infinitely extended fluid reservoir.
All discussions of the previous sections have been based on the assumption that all
subsets of $\{x_3 \leq 0\}\setminus \Omega^y$ are actually filled with
fluid, see Figure~\ref{figure1}. This is nonrestrictive in case
$\Omega^y$ is convex (a stable property with respect to weak $W^{1,p}$
convergence of deformations). Still, the case of a nonconvex
$\Omega^y$ is of  major  applicative relevance, for it corresponds to
the idealized situation of a floating vessel, see Figure~\ref{ship}.
In this
case, the set $\{x_3 \leq 0\}\setminus \Omega^y$ may have different
connected components, one of which is unbounded (since $\Omega^y$ is
necessarily bounded). We shall use the notation
$$D^y:=\text{union of all bounded connected components of $\{x_3 \leq
  0\}\setminus \ove {\Omega^y}$}.$$ 
The former choice in \eqref{eq:E} corresponds to  the assumption   that $D^y$ is
filled by the fluid. On the other hand, one could assume that some of
the connected components of $D^y$ do not contain fluid, or are
partially filled with fluid. Among  all  options, we consider
here the case in which $D^y$ contains no fluid at all, see Figure~\ref{ship}. In this setting, the energy $E$ is redefined as 
\begin{equation}\label{eq:E2}
\RRR E_{\rm ship}(y) \EEE = \int_\Omega W(\nabla y(X))\, \dX +
  \int_{\Omega^y\cup D^y} g\rho_f x_3^-\, \dx +
\int_\Omega g \rho_sy_3(X)\, \dX.
\end{equation}
In particular, the set $D^y$ contributes to the hydrostatic term of
the energy, for the fluid is displaced out of $D^y$, although $D^y$ is
not a subset of $\Omega^y$.  For a justification in terms of the
Euler-Lagrange equations we refer to Section \ref{sec:basic}, 
where now  \eqref{eq:equilibrium4}  holds  with $\partial (\Omega^y\cup
D^y)$  in place of $\partial \Omega^y$.  

The effect  of
keeping track of the set $D^y$ in the integral of the hydrostatic term
allows    the
solid to float, even for $\rho_s>\rho_f$, as  it  commonly happens for
usual vessels.  From now on, we restrict our considerations to the case $\rho_s>\rho_f$ and suppose that the solid is incompressible, i.e., $J^y=1$ almost everywhere.  Note that then the energy is again unbounded
from below since the energy decreases as the solid sinks, as soon as
$ \sup  y_3 <0$ (and hence $D^y=\emptyset$). This is indeed the
mechanism responsible for all shipwrecks.

As no global minimizers can be expected to exist, we address the
existence of floating  solids locally
minimizing  the energy \RRR $E_{\rm ship}$ \EEE \eqref{eq:E2}.  Local minimizers can only be expected if $|D^y|$ is sufficiently large. Therefore, we impose 
\begin{equation}
\label{eq:con}
|D^y| \geq \eta \quad\text{for} \quad \eta = \frac{\rho_s - \rho_f}{\rho_f}|\Omega|.
\end{equation}
The specific choice of $\eta$ is  indeed  tailored to let   $|D^y|=\eta$
exactly correspond to the case of barely floating solids.  In
fact,  one can extend the discussion leading to \eqref{eq:archi}
to the specific case of energy  \RRR $E_{\rm ship}$ \EEE from  \eqref{eq:E2} in order to  derive   that a
critical point  of \RRR $E_{\rm ship}$  \EEE  fulfills
$$g\rho_f|(\Omega^y \cap \{x_3 \leq 0 \} )\cup D^y| = g\rho_f
 |\Omega^y \cap \{x_3 \leq 0 \}| + g\rho_f |D^y|  = g\rho_s|\Omega|.$$
In case $\Omega^y$ is barely floating, since $|\Omega^y|=|\Omega|$, one
has that 
$$ \rho_f |\Omega| +\rho_f |D^y| = \rho_s |\Omega|$$
and  therefore  $|D^y|=\eta$. On the other hand, if $|D^y|=\eta$, we get 
$$ |\Omega^y \cap \{x_3 \leq 0 \}|=
(\rho_s/\rho_f)|\Omega|  - \eta  = |\Omega|=|\Omega^y|,$$
which implies that $\Omega^y$ is barely floating.  We  have
hence proved  that  in
case of $|D^y|>\eta$ the solid necessarily floats.

\begin{figure}[h]
  \centering
  \pgfdeclareimage[width=120mm]{ship2}{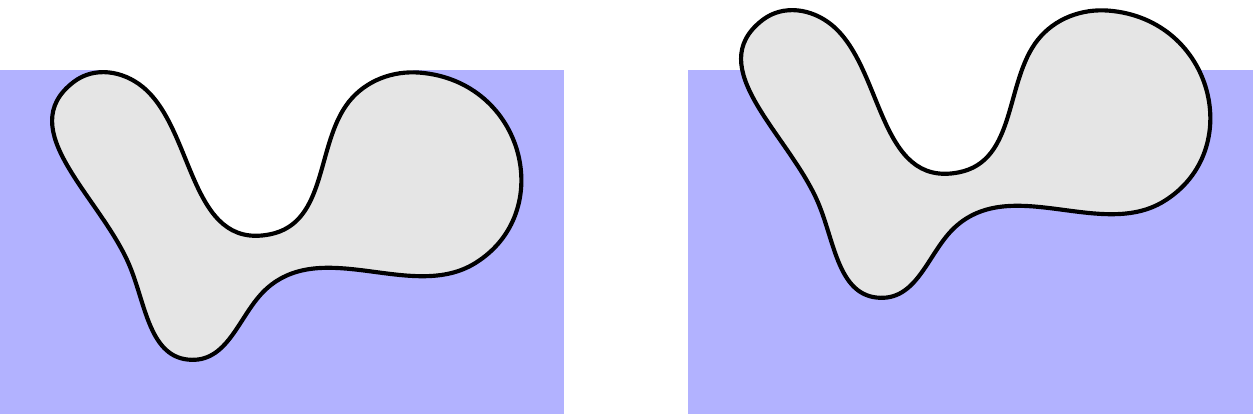}
    \pgfuseimage{ship2}
  \caption{A barely floating solid (left) and an admissible $y\in A$
    with $A$ from \eqref{eq:A} (right).}
\label{ship2}
\end{figure}

As in the previous sections, our goal is to show the existence of (local) minimizers by the \BBB direct method. \EEE  The presence of the
set $D^y$ in \eqref{eq:E2}, however,  is posing lower semicontinuity
problems: let  $\Omega^y$  be the barely floating deformation depicted on the
left of Figure \ref{ship2}
and consider the sequence $y_k = y - e_3/k$. As $|D^y|>0$ but
$|D^{y_k}|=0$, one has that 
\begin{align*}&\liminf_{k\to \infty} \int_{\Omega^{y_k}\cup D^{y_k}}g \rho_f x_3^-
\, \dx = \liminf_{k\to \infty} \int_{\Omega^{y_k}}g \rho_f x_3^-
\, \dx  \\
&\quad =  \int_{\Omega^{y}}g \rho_f x_3^- \, \dx
<\int_{\Omega^{y}\cup D^{y}}g \rho_f x_3^- \, \dx 
\end{align*}
and lower semicontinuity fails. 
 The case of a barely floating solid is, however, of  limited practical interest. 

In order to restore  lower semicontinuity, we strengthen the
requirements by seeking local minimizers such that  the  map  $
y\mapsto |D^y|$  is continuous in a specific way, namely  we ask
 $y$  to satisfy  for all $\varepsilon \in (0,\varepsilon_0)$ ($\varepsilon_0>0$ fixed)
\begin{align}\label{eq: new condition}
  \bar{y}\in W^{1,p}(\Omega;\Rz^3)\colon \quad   \| y_3 - \ove
  y_3\|_{L^\infty(\Omega)}<\epsi \ \Rightarrow \ |D^{\ove y} \triangle
  D^y| \le \varepsilon^{1/2}. 
\end{align} 
We now define the set of admissible deformations as 
\begin{equation}A=\{y \in W^{1,p}(\Omega;\Rz^3)\ \ \text{fulfilling \eqref{eq:con} and \eqref{eq: new condition}}\}.
\label{eq:A}
\end{equation}
 Note that $A$ is not empty, for any reference configuration $\Omega$ can
be deformed into a thin spherical  half-shell  of arbitrarily given internal radius,
possibly at a large elastic energy cost.  In this case, \eqref{eq:con} and \eqref{eq: new condition} can be indeed verified.  

As all solids $\Omega^y$ with $y\in A$ are   floating or barely floating,  i.e.,  $\sup y_3\ge 0$,   the energy \RRR $E_{\rm ship}$ \EEE
  is coercive on $A$  by \eqref{eq:diam}--\eqref{eq:upper}.  Moreover, it is elementary to check that $A$ is closed under uniform convergence.    In order to check lower
semicontinuity, we see that we can  pass to the limit in the hydrostatic terms by \eqref{eq: new condition}. This yields the existence of  a minimizer in $A$.

Let  $y^*\in A$ be a minimizer of \RRR $E_{\rm ship}$ \EEE in $A$.  We have that 
$\Omega^{y^*}$ 
is floating but not barely floating, see Figure~\ref{ship2}
right. Indeed, assume  
$\Omega^{y^*}$  to be barely floating. Then, $\ove y  =
y^*-(\varepsilon_0/2)e_3$ would have $D^{\ove y} = \emptyset$ and
$|D^{\ove y} \triangle D^{y^*}|=|D^{y^*}|\geq \eta$ would entail that 
$y^*$ does not satisfy \eqref{eq: new condition}, provided $\varepsilon_0$ is small with respect to $\eta$. This is a contradiction, and we  have hence proved that
$|D^{y^*}|>\eta$. 

 Moving from this, one has that $y^*$ is a true local
minimizer of \RRR $E_{\rm ship}$ \EEE if there exists  $\epsi_0>0$ such that
for all $y$ with $\| y^*_3 - y_3\|_{L^\infty(\Omega)}<\epsi<\epsi_0$ one has
that $|D^y \triangle D^{y^*}|<\epsi^{1/2}$. An example in this
direction is in Figure~\ref{ship2} right.

\begin{figure}[h]
  \centering
  \pgfdeclareimage[width=60mm]{ship4}{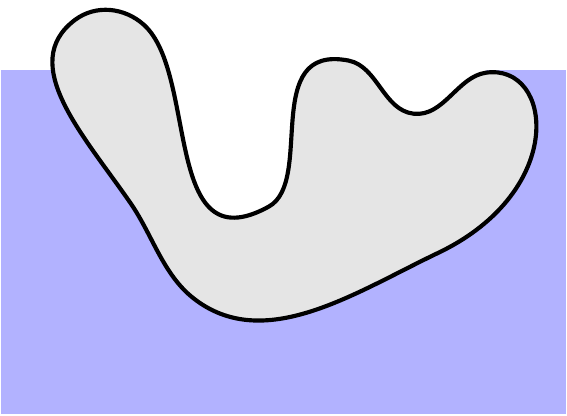}
    \pgfuseimage{ship4}
  \caption{A deformation with $\omega^{y^*}\subset\Omega^{y^*}$ with
$\sup_{\omega^{y^*}} y^*_3 =0$ and  $\sup_{\Omega^{y*}}y^*_3>0$}
\label{ship4}
\end{figure}

 Note that the case of $y^*$
not being a local minimizer, for all $\epsi_0>0$, is of no real
applicative interest. In fact, this would happen if one could find a
sequence $\epsi_n \to 0$ and deformations $y^n$ such that $\| y^*_3 -
 y^n_3\|_{L^\infty(\Omega)}<\epsi_n$ but $|D^{y^n} \triangle
D^{y^*}|\geq \epsi_n^{1/2}$.   We see, in this case,  that small changes in the deformation  cause large variations  in  the volume of the  symmetric difference and this     means that $\Omega^{y^*}$  is either barely floating, see Figure~\ref{ship2} left, or that
a part of the ship is barely floating, i.e., the picture applies only
to an open subset $\omega^{y^*}\subset\Omega^{y^*}$ for which
$\sup_{\omega^{y^*}} y^*_3 =0$ although  $\sup_{\Omega^{y*}}y^*_3>0$,
see Figure~\ref{ship4}. 


Let us conclude by explicitly remarking that we are not in the
position of establishing a priori if  $y^*$ is indeed a local
minimizer (for some  $\epsi_0 >0$)  or not. Indeed, this depends on
the data of the problem, most notably on $\Omega$ and $W$. A positive
example would be given by a   spherical half-shell  $\Omega$ with
sufficiently large radius under an energy density $W$ very much penalizing deformations far from
identity. On the other hand, the same $W$ in case of $\Omega$
being a ball can be expected to provide a negative example.

\UUU
\section{Outlook}\label{sec:outlook}

We conclude our discussion by mentioning some possible future
developments.
A number of interesting extensions are indeed within easy reach. One could for
instance consider the case of non homogeneous densities $\rho_s
\colon \Omega \to (0,\infty)$ and $\RRR \rho_f \colon \BBB \Rz^3 \to (0,\infty)$ (this
second one possibly depending on $x_3$ only). Also the case of a
nonhomogeneous gravity $g\colon \Rz^3 \to (0,\infty)$ could be easily
handled. 

For a porous solid, the analysis can be readily extended to the
capillarity case, if the height of the wet part of the body over \RRR water level \BBB is a priori fixed. This can be expected to be quite common for relatively large
systems under moderate deformations (a wooden raft, for instance). For smaller
systems or large strains (a soft sponge), capillarity depends on the
deformation, adding a challenging coupling effect to the picture. The interplay between
capillarity and anisotropy may also
enter the picture, making the analysis even more involved.

Still in the stationary regime (evolution is out of the scope of the paper), one could
consider the fluid to be viscous, incompressible, and in stationary
motion. This amounts to couple the minimization of the energy with the
stationary Navier-Stokes or with the Stokes system. In this case, the hydrostatic
pressure term has to be augmented by the normal component of the
normal stress. To the best of our knowledge, an existence result in
this direction is still missing. The reader is however referred to \cite{Seb}, where the dynamic setting is tackled.

The equilibrium problem for the immersed body can be combined with
other physical effects, in a multiphysics setting. Thermal and
electromagnetic effects can be taken into the picture, for
instance.
A particularly relevant prospect would be to allow for contact with
other hyperleastic bodies or with the bottom of the reservoir. Here,
one would be asked to limit the set of admissible configurations to
those avoiding material interpenetration, an option which should be
amenable, possibly under stricter conditions \cite{Pantz}. 

Another relevant issue is the stability of local minimizers  (imagine a loaded ship at sea). This issue is classical in the
case of a rigid convex body \cite{Guerrero,Kurusa,Wegner} and, to our knowledge, completely open in
the case of a deformable body. Note once again that our analysis of
the ship situation from Section \ref{sec:ship} left open the question
on how to  check local minimality of the configuration $y^*$, given 
information on $\Omega$ and $W$.

Finally, topology optimization problems can also be considered for partly or fully immersed bodies. Here we have in mind e.g.~the  design of completely submerged buoyant modules with design-dependent fluid pressure loading. This type of structure is used to
support offshore rig installation and pipeline transportation at various  water
depths. The  optimization methods seek to identify the buoy design with  the highest stiffness and  allowing  to withstand deepwater
pressure. At the same time,  it should use  the least material and \RRR have \BBB a minimum prescribed buoyancy. We refer to \cite{Picelli} for  details\EEE.

\section*{Acknowledgments}
 MF was supported by  the  DFG-FWF project FR 4083/3-1/I\,4354,  and by the Deutsche Forschungsgemeinschaft (DFG, German
  Research Foundation) under Germany's Excellence Strategy EXC
  2044-390685587, Mathematics M\"unster:
  Dynamics--Geometry--Structure. 
 US was   supported by the Vienna Science and Technology Fund
(WWTF) through the project MA14-009 and by the Austrian Science Fund
(FWF) project F\,65.  MK was supported by the GA\v{C}R-FWF project 
19-29646L-I\,4052.


\begin{thebibliography}{99}


\bibitem{Archi}
T. L. Heath (Ed.). {\it The works of Archimedes}. Cambridge
University Press, 1897. Reprinted Dover, Mineola, NY, 2002. 

\BBB
\bibitem{Seb}
B. Bene\v sov\'a, M. Kampschulte, S. Schwarzacher.
A variational approach to hyperbolic evolutions and fluid-structure
interactions. \href{https://arxiv.org/abs/2008.04796}{\tt arxiv:2008.04796}.

\EEE

\bibitem{Bishop}
R.~E.~D.~Bishop, W.~G.~Price. {\it Hydroelasticity of ships}. Cambridge University Press, 1979.


\bibitem{CN}
P. G. Ciarlet, J. Ne\v cas.
Injectivity and self-contact in nonlinear elasticity. {\it Arch. Ration. Mech. Anal.} 97 (1987), 171--188. 



\bibitem{Finn0}

R. Finn. Floating and partly immersed balls in a weightless
environment. {\it Funct.
Differ. Equ.} 12 (2005), 167--173.

\bibitem{Finn}
R. Finn. Criteria for Floating I.
{\it 
J. Math. Fluid Mech.} 13 (2011), 103--115.

\bibitem{Finn2}
R. Finn, T.~I. Vogel. Floating criteria in three
dimensions. {\it Analysis (Munich)}, 29 (2009), 387--402.
Erratum. {\it Analysis (Munich)}, 29 (2009), 339. 



\bibitem{GKMS}
 D. Grandi, M. Kru\v z\'\i k, E. Mainini, U. Stefanelli. A phase-field
 approach to interfacial energies in the deformed configuration. {\it
   Arch. Ration. Mech. Anal.} 234 (2019), 351--373.

\bibitem{Guerrero}
Z. Guerrero-Zarazua, J. Jer\'onimo-Castro. Some comments on floating
and centroid bodies in the plane. {\it Aequationes Math.} 92 (2018), 211--222. 

\bibitem{HK}
 S. Hencl, P.   Koskela. {\it Lectures on mappings of finite
   distortion}. Lecture Notes in Mathematics 2096, Springer, 2014.


\bibitem{compos2}
 M. Kru\v z\'\i k, D. Melching, U. Stefanelli. {Quasistatic
   evolution for dislocation-free finite plasticity}. {\it ESAIM Calc. Var. Control Optim.} DOI: https://doi.org/10.1051/cocv/2020031  \href{https://arxiv.org/abs/1912.10118}{\tt arXiv:1912.10118}.

\bibitem{John}
F. John. On the motion of floating bodies, I, II. {\it Comm. Pure Appl. Math.} 2
(1949), 13--57 \& 3 (1950), 45--101.

\bibitem{Kaltenbacher}
B. Kaltenbacher, I. Kukavica,  I. Lasiecka, R. Triggiani, A. Tuffaha,
J.~T. Webster. {\it Mathematical theory of evolutionary fluid-flow structure interactions}. Lecture notes from Oberwolfach seminars, November 20–26, 2016. Oberwolfach Seminars, 48. Birkhäuser/Springer, Cham, 2018.

\bibitem{Kruzik-Roubicek}
M. Kru\v z\'\i k, T.~Roub{\'i}{\v c}ek. {\it Mathematical methods in
  continuum mechanics of solids}. Interaction of Mechanics and
Mathematics. Springer, Cham, 2019.

\bibitem{Kurusa}
A. Kurusa,T. \'Odor. Spherical floating bodies. {\it Acta
  Sci. Math. (Szeged)}, 81 (2015),   699--714.

\bibitem{Laplace}
P.~S. Laplace. {Trait\'e de m\'ecanique c\'eleste: supplement 2,
  909--945, au Livre X}. In {\it Oeuvres Compl\`ete}, vol. 4. Gauthier
Villars, Paris. English translation by N. Bowditch (1839), reprinted
by Chelsea, New York, 1966.

\bibitem{Mauldin}
  R.~D. Mauldin (ed.). {\it The Scottish Book}. Birkh\"auser, Boston, 1981.



\bibitem{McCuan}
J. McCuan. {A variational formula for floating bodies}.
 {\it Pacific J. Math.}  231 (2007), 167--191.


\bibitem{McCuan2}
J. McCuan.
Archimedes Revisited. {\it 
Milan J. Math.} 77 (2009), 385--396.

\bibitem{McCuan3}
J. McCuan, R. Treinen. Capillarity and Archimedes' principle of
flotation. {\it Pacific J. Math.} 265 (2013),  123--150.

\BBB
\bibitem{Pantz}
O. Pantz. The modeling of deformable bodies with frictionless
(self-)contacts. {\it Arch. Ration. Mech. Anal.} 188 (2008),   183--212.
\EEE

\BBB
\bibitem{Picelli}
R. Picelli, R.van Dijk, W.M. Vincente, R. Pavanelloa, M. Langelaar, F. van Keulen. Topology optimization for submerged buoyant structures. {\it  Engrg. Optim.} 49 (2017), 1--21.
\EEE


\bibitem{Richter}
T. Richter. {\it  Fluid-structure interactions. Models, analysis and
  finite elements}. Lecture Notes in Computational Science and
Engineering, 118. Springer, Cham, 2017.

\bibitem{Treinen}
R. Treinen. A general existence theorem for symmetric floating
drops. {\it Arch. Math. (Basel)}, 94 (2010),  477--488.

\bibitem{Wegner}
F. Wegner.
Floating bodies of equilibrium. 
{\it Stud. Appl. Math.} 111 (2003),  167--183.

\end{thebibliography}
\end{document}